\documentclass{proc-l}
\usepackage{amsmath}
\usepackage{amsfonts}
\usepackage{amssymb}

\newcommand\eps{\varepsilon}
\newcommand\R{{\mathbf{R}}}

\newcommand\Z{{\mathbf{Z}}}

\renewcommand\S{{\mathcal{S}}}
\newcommand\E{{\mathcal{E}}}
\newcommand\M{{\mathcal{M}}}
\renewcommand\L{{\mathcal{L}}}

\newtheorem{theorem}{Theorem}[section]
  \newtheorem{conjecture}[subsection]{Conjecture}
  \newtheorem{proposition}[subsection]{Proposition}
  \newtheorem{lemma}[subsection]{Lemma}
  \newtheorem{corollary}[subsection]{Corollary}

\theoremstyle{definition}

\theoremstyle{remark}
\newtheorem{remark}[theorem]{Remark}

\numberwithin{equation}{section}

\begin{document}

\title[Cone restriction]{A note on the cone restriction conjecture in the cylindrically symmetric case}
\author{Shuanglin Shao}
\address{Department of Mathematics, UCLA, Los Angeles CA 90095-1555}
\email{slshao@math.ucla.edu}


\date{October 8, 2007 and, in revised form, June 10, 2008.}

\commby{Hart F. Smith}

\begin{abstract}
In this note, we present two arguments showing that the classical \textit{linear adjoint cone
restriction conjecture} holds for the class of functions supported on the cone and invariant under the spatial rotation in all dimensions. The first is based on a dyadic restriction estimate, while the second follows from a strengthening version of the Hausdorff-Young inequality and the H\"older inequality in the Lorentz spaces.
\end{abstract}

\maketitle

\section{Introduction}
Let $n\ge 2$ be a fixed integer and $S$ be a smooth compact non-empty subset of the cone
$\{(\tau,\xi)\in \R\times \R^n:\,\tau=|\xi|\}$, where we interpret $\R\times \R^n$ as the time-space
frequency space. If $0<p,q\le \infty$, the classical \emph{linear adjoint restriction estimate}
\footnote{In the notation of \cite{Tao:2008:recent-progress-restric}, the estimate \eqref{eq:lin-restr}
 is denoted by $R^{*}_S(p\to q)$.} for the cone is the following \textit{``a priori"} estimate
\begin{equation}\label{eq:lin-restr}
\|(fd\sigma)^{\vee} \|_{L_{t,x}^q(\R\times \R^n)}\le C_{p,q,n, S}\|f\|_{L^p(S,d\sigma)}
\end{equation} for all Schwartz functions $f$ on $S$, where
$$(fd\sigma)^{\vee}(t,x)=\int_{S} f(\tau, \xi)e^{i(x\cdot\xi+t\tau)}d\sigma(\xi)
=\int_{\R^n}f(|\xi|, \xi)e^{i(x\cdot\xi+t|\xi|)}\frac {d\xi}{|\xi|}$$ denotes the inverse space-time
Fourier transform of the measure $fd\sigma$, and $d\sigma$ is the pull-back of the measure $\frac {d\xi}{|\xi|}$ under the projection map $(\tau,\xi)\mapsto \xi$. By duality, the estimate
\eqref{eq:lin-restr} is equivalent to
\begin{equation*}
\|\hat{f}\|_{L^{p'}(S, d\sigma)}\le C_{p,q,n,S} \|f\|_{L^{q'}(\R\times \R^n)}
\end{equation*} for all Schwartz functions $f$, which roughly says that the Fourier transform
of an $L^{q'}(\R\times \R^n)$ function can be ``meaningfully" restricted to the cone
$S$. This leads to the \textit{restriction problem}, one of the central problems in harmonic
analysis, which concerns the optimal range of exponents $p$ and $q$ for which the estimate
\eqref{eq:lin-restr} should hold. It was originally proposed by Stein for the sphere
\cite{Stein:1979:problems-in-harmonic} and then extended to smooth sub-manifolds of $\R\times
\R^n$ with appropriate curvature \cite[Chapter 8, pages 352-355, 364-367 ]{Stein:1993} such as the paraboloid and the cone. The restriction problem is intricately related to other outstanding problems in analysis such as the Bochner-Riesz conjecture, the Sogge's local smoothing conjecture, the Kakeya set conjecture
and the Kakeya maximal function conjecture, see e.g., \cite{Tao:2008:recent-progress-restric},
\cite{Tao:1999:Bochner-Resiez-restri}.

By testing \eqref{eq:lin-restr} against the characteristic functions supported on a symmetric band or a small cap of the cone, the following conjecture on the restriction of the Fourier transform to the cone can be formulated,
\begin{conjecture}[Linear adjoint cone restriction conjecture]\label{con:lin-restr}
The inequality \eqref{eq:lin-restr} holds with constants depending on $S$, $n$ and $p,q$ if and
only if $q>\frac {2n}{n-1}$ and $\frac {n+1}{q}\le \frac {n-1}{p'}$.
\end{conjecture}

C\'ordoba and Stein proved that \eqref{eq:lin-restr} was true under the condition $p=2$ and $q\ge \frac {2(n+1)}{n-1}$ in an unpublished work. Strichartz \cite{Strichartz:1977} then extended the results to more general quadratic surfaces. In 1985, Barcelo \cite{Barcelo:1985:restic-cone} proved Conjecture \ref{con:lin-restr} when $n=2$. A major breakthrough was made in 2001 by Wolff \cite{Wolff:2001:restric-cone}, who showed that Conjecture \ref{con:lin-restr} was true when $n=3$. This was based on a new bilinear cone restriction estimate, which also gave the current best result $q> \frac {2(n+3)}{n+1}$ in higher dimensions $n\ge 4$. We should remark that all the recent progress on the linear restriction is achieved from the corresponding bilinear restriction estimates, especially the bilinear $L^2$-type estimates, $L^2\times L^2\to L^{q}$ for some $q\in [1,2]$; more information about the so-called bilinear method and recent ideas of attacking the restriction conjecture such as the reduction to the local restriction estimates, the wave packet decomposition and the induction-on-scales can be found in \cite{Tao-Vargas-Vega:1998:bilinear-restri-kakeya}, \cite{Tao:2001:endpoint-cone}, \cite{Wolff:2001:restric-cone}, \cite{Tao:2003:paraboloid-restri} and \cite{Tao:2008:recent-progress-restric}.

When we restrict the functions supported on the cone $S$ to be cylindrically symmetric, i.e., functions invariant under the spatial rotation, the following theorem is our main result in this paper,
\begin{theorem}\label{thm:line-con} Conjecture \ref{con:lin-restr}
holds for cylindrically symmetric functions supported on the cone in all dimensions.
\end{theorem}

For the same class of functions but supported on the paraboloid, the author \cite{Shao:2008:paraboloid} has verified the corresponding conjecture for the paraboloid in all dimensions. The first proof of Theorem \ref{thm:line-con} is along similar lines as in \cite{Shao:2008:paraboloid}, through dyadically decomposing both the frequency and spatial spaces and then establishing a family of dyadical restriction estimates based on the ``Fourier-Bessel" formula defined in Section \ref{sec:dya-pf}. The second proof is inspired by Nicola's argument on the implication of cone restriction conjecture from the sphere restriction conjecture in \cite{Nicola:2008-cone-sphere-restriction}. The key ingredient is the use of the strengthening version of the Hausdorff-Young inequality \cite[Chapter 4, Corollary 3.16]{Stein-Weiss:1971:fourier-analysis} and the H\"older inequality in the Lorentz spaces \cite[Chapter 5, Theorem 5.3.1]{Bergh-Lofstrom:1976-interpolation-space}.

\begin{remark}\label{re:more-dydic-estmte}
As in \cite{Shao:2008:paraboloid}, for the cylindrically symmetric functions with dydical supports,
we expect that more estimates are available. This is indeed the case:  from Proposition
\ref{prop:sml-R} and Corollary \ref{cor:big-R}, when $f$ is cylindrically symmetric and supported on a subset of the cone $\{(|\xi|,\xi):1\le |\xi|\le 2\}$,  for $q>\frac{2n}{n-1}$ and $ q\ge p'$,
\begin{equation*}
\|(fd\sigma)^\vee\|_{L^q(\R\times \R^n)} \le C_{p,q} \|f\|_{L^p(S,d\sigma)}.
\end{equation*}
We note that $q\ge p'$ is an improvement over $q\ge \frac{n+1}{n-1}p'$.
\end{remark}

\begin{remark}
When $S$ is the whole cone instead of a compact subset of the cone, we see that the necessary conditions are strengthened to
$$q>\frac {2n}{n-1},\quad \frac {n+1}{q}=\frac{n-1}{p'}.$$
In this case, on the one hand, Theorem \ref{thm:line-con} guarantees that the cone restriction conjecture \ref{con:lin-restr} is true; on the other hand, unlike the situation in Remark \ref{re:more-dydic-estmte}, there are no more estimates available.
\end{remark}

This paper is organized as follows. Section \ref{sec:nota} is devoted to establishing the standard
notations; in Section \ref{sec:dya-pf} we present our first proof of Theorem \ref{thm:line-con} via
the dyadic restriction estimates; in Section \ref{sec:Nic'pf} we present another proof by using a strengthening version of the Hausdorff-Young inequality and the H\"older inequality in the Lorentz spaces.

$\textbf{Acknowledgments.}$ The author is very grateful to his advisor, Terence Tao, for the helpful discussions on this problem, for his support during the preparation of this paper. The author would also like to thank the referee for their valuable suggestions and comments.

\section{Notations}\label{sec:nota}
We will use the notations $X\lesssim Y$, $Y\gtrsim X$, or $X=O(Y)$ to denote the estimate $|X|\le
C Y$ for some constant $0<C<\infty$, which may depend on $p,q,n$ and $S$, but not on the
functions. If $X\lesssim Y$ and $Y\lesssim X$ we will write $X\sim Y$. If the constant $C$ depends
on a special parameter other than the above, we shall denote it explicitly by subscripts. For
example, $C_{\eps}$ should be understood as a positive constant not only depending on $p,q,n$ and
$S$, but also on $\eps$.

By $\S^{n-1}$ we denote the ${n-1}$ dimensional unit sphere, and by $d\mu$ the canonical surface measure of the sphere. We define a dyadic number to be any number $R\in 2^{\Z}$ of the form $R=2^j$ where $j$ is an integer. For each dyadic number $R>0$, we define the dyadic annulus in $\R^{n}$, $A_R:=\{x\in \R^{n}: R/2\le |x|\le R\}$. By $\L_N$, we denote the class of cylindrically symmetric functions dyadically supported on the cone, i.e., functions invariant under the spatial rotation and supported on a set of the form $\{(\tau,\xi): N\le |\xi|\le 2N, \tau=|\xi|\}$ with dyadic $N>0$.  We define the spacetime norm $L^q_tL^r_x$ of $f$ on $\R\times \R^{n}$ by $$\|f\|_{L^q_tL^r_x(\R\times\R^{n})}:=\left(\int_{\R}\left(\int_{\R^{n}}|f(t,x)|^{r}d\,x\right)^{q/r}
dt\right)^{1/q}$$ with the usual modifications when $q$ or $r$ are equal to infinity, or when the domain $\R\times \R^{n}$ is replaced by a small region of spacetime such as $\R\times A_R$. When $q=r$, we
abbreviate it by $L^{q}_{t, x}$.  We define the spatial Fourier transform of $f$ on $\R^{n}$ by $\hat{f}(\xi)=\int_{\R^{n}}f(x)e^{-ix\cdot\xi}dx$. We use $1_U$ to denote the characteristic function of the set $U$, i.e., $1_U(x):=1$ if $x\in U$, otherwise $0$. For $1\le p\le \infty$, we denote the conjugate exponent of $p$ by $p'$, i.e.,  $1/p+1/p'=1$.

\section{First proof of Theorem \ref{thm:line-con}}\label{sec:dya-pf}
For any cylindrically symmetric function $f$ on the cone, we set $F(|\xi|):=f(|\xi|, \xi)$. We
observe that $(fd\sigma)^\vee(t, x)$ is also a cylindrically symmetric function. To begin the
proof of Theorem \ref{thm:line-con}, we investigate the behavior of $(fd\sigma)^\vee$ on
$\{|x|\le 1\}$ via the following proposition.
\begin{proposition}\label{prop:sml-R}
Suppose $f\in \L_1$. Then for any $1\le p\le \infty$, $q\ge \max\{2, p'\}$ and $R\le 1$, we have
\begin{equation}\label{eq:err-1}
\|(fd\sigma)^\vee\|_{L^q_{t,x}(\R\times A_R)} \lesssim R^{\frac{n}{q}}\|f\|_{L^p(S,d\sigma)}.
\end{equation}
\end{proposition}

\begin{proof}
If we change to polar coordinates, the left-hand side of \eqref{eq:err-1} is
\begin{align*}
&\quad\left(\int_{A_R}\int_{\R} \left|\int_{1 \le|\xi|\le 2}f(|\xi|,\xi) e^{i(x\cdot\xi+t|\xi|)}\frac{d\xi}{|\xi|}\right|^q dtdx \right)^{1/q}\\
&=\left(\int_{R/2}^{R}\int_{\R}\left|\int_{1 \le|\xi|\le 2}f(|\xi|,\xi) e^{i(re_1\cdot\xi+t|\xi|)}\frac{d\xi}{|\xi|}\right|^q dt\,r^{n-1}dr \right)^{1/q}\\
&=\left(\int_{R/2}^{R}\int_{\R}\left|\int_{I}
F(s)s^{n-2}e^{its}\int_{\S^{n-1}}e^{irse_1\cdot\omega}d\mu(\omega)ds\right|^q dt\,r^{n-1}dr \right)^{1/q}\\
&=\left(\int_{R/2}^{R}\int_{\R}\left|\int_{I}
F(s)s^{n-2}e^{its}(d\mu)^\vee(rse_1\cdot\omega) ds\right|^q dt\,r^{n-1}dr \right)^{1/q},
\end{align*}
where $I=[1,2]$, $e_1=(1,0,\ldots,0)\in \R^n$ and $``\cdot"$ denotes the inner product operation in $\R^n$. Then from the Hausdorff-Young inequality when $q>2$ or Plancherel theorem when $q=2$ and using $\|(d\mu)^\vee\|_{L^{\infty}_{\omega}}\lesssim 1$, the left-hand side of \eqref{eq:err-1} is further bounded by
\begin{equation*}
R^{\frac{n-1}{q}}\left(\int_{R/2}^R \|F\|^q_{L^{q'}(I)}\,dr\right)^{1/q}\sim
R^{\frac{n}{q}}\|F\|_{L^{q'}(I)}.
\end{equation*}
Then by applying the H\"older inequality to raising $q'$ to $p$ since $p\ge q'$, and noting $\|F\|_{L^{p}(I)}\sim \|f\|_{L^p(S,d\sigma)}$,
\eqref{eq:err-1} follows.
\end{proof}

Before investigating the behavior of $(fd\sigma)^\vee$ on $|x|\ge 1$, we exploit the
cylindrical symmetry of $f$ in the following proposition. Note that we will encode the error term of the Bessel function into integrals instead of using its asymptotic bound.
\begin{lemma}[Fourier-Bessel formula]\label{le:fouri-bessl}
Suppose $f$ is a cylindrically symmetric function supported on the cone. Then
\begin{align*}\label{eq:Fourier-Bessel}
&(fd\sigma)^\vee(t,x)\\
&=c_nr^{-\frac{n-1}{2}}\int_{I} F(s)s^{\frac{n-3}{2}} e^{i(rs+ts)}ds
+c_nr^{-\frac{n-1}{2}}\int_{I} F(s)s^{\frac{n-3}{2}} e^{i(-rs+ts)}ds\\
&\quad+c_n\int_{I}F(s)s^{n-2}e^{its-irs}\int_{0}^{\infty}e^{-rsy}y^{\frac{n-3}{2}}
[(y+2i)^{\frac{n-3}{2}}-(2i)^{\frac{n-3}{2}}]dyds\\
&\quad+c_n\int_{I}F(s)s^{n-2}e^{its+irs}\int_{0}^{\infty}e^{-rsy}y^{\frac{n-3}{2}}
[(y-2i)^{\frac{n-3}{2}}-(-2i)^{\frac{n-3}{2}}]dyds.
\end{align*}
where $I$ denotes the interval in the radial direction and $r=|x|$.
\end{lemma}
\begin{proof}
We first expand $(fd\sigma)^\vee$ in the polar coordinates,
\begin{equation*}
(fd\sigma)^\vee(t, x)=\int_{\{|\xi|\in I\}} f(|\xi|,
\xi)e^{i(re_1\cdot\xi+t|\xi|)}\frac{d\xi}{|\xi|}=\int_I F(s)e^{its}s^{n-2}(d\mu)^\vee(rse_1\cdot\omega)ds.
\end{equation*}
We recall $(d\mu)^\vee(\xi)=c_n|\xi|^{\frac {2-n}{2}}J_{\frac{n-2}{2}}(|\xi|),$
see e.g., \cite[page 347]{Stein:1993}. Moreover from \cite[Chapter 3, Lemma 11]{Stein-Weiss:1971:fourier-analysis}, we obtain that, for fixed $m\ge 0$,
\begin{align*}
J_m(r)&=c_mr^{-1/2}(e^{ir}-e^{-ir}) \\ &+c_mr^{m}e^{-ir}\int_{0}^{\infty}e^{-ry}y^{\frac{2m-1}{2}}[(y+2i)^{\frac{2m-1}{2}}
-(2i)^{\frac{2m-1}{2}}]dy\\
&+c_mr^me^{ir}\int_{0}^{\infty}e^{-ry}y^{\frac{2m-1}{2}}[(y-2i)^{\frac{2m-1}{2}}
-(-2i)^{\frac{2m-1}{2}}]dy.
\end{align*}
Then Lemma \ref{le:fouri-bessl} follows after we combine these two estimates and
set $m=\frac{n-2}{2}$.
\end{proof}

In view of the previous lemma, we thus define the main term and the error term of $(fd\sigma)^\vee$ by
\begin{align*}
\mathcal{M}f(t,x)&:=c_nr^{-\frac{n-1}{2}}\int_{I}F(s)s^{\frac{n-3}{2}} e^{i(rs+ts)}ds+
c_nr^{-\frac{n-1}{2}}\int_{I}F(s)s^{\frac{n-3}{2}} e^{i(-rs+ts)}ds,\\
\mathcal{E}f(t,x)&:=c_n\int_{I}F(s)s^{n-2}e^{its+irs}\int_{0}^{\infty}e^{-rsy}y^{\frac{n-3}{2}}
[(y+2i)^{\frac{n-3}{2}}-(2i)^{\frac{n-3}{2}}]dyds\\
&-c_n\int_{I}F(s)s^{n-2}e^{its-irs}\int_{0}^{\infty}e^{-rsy}y^{\frac{n-3}{2}}[(y-2i)^{\frac{n-3}{2}}
-(-2i)^{\frac{n-3}{2}}]dyds.
\end{align*}
Heuristically, one should think of $\E f$ as $r^{-(n+1)/2}\int_{I}F(s)s^{\frac{n-5}{2}}e^{its}ds$, which is given by estimating the error term of Bessel function $J_m(r)$ by $r^{-3/2}$. The following
proposition shows that the error term estimate is acceptable compared to the main term estimate.
\begin{proposition}\label{prop:big-R}
Suppose $f\in \L_1$. Then for all $1\le p\le \infty$, $q\ge \max\{2, p'\}$, a dyadic number $R\ge 2$ and $f\in L^p(S,d\sigma)$, we have the main term estimate,
\begin{equation}\label{eq:main-estmt}
\|\M f\|_{L^q_{t,x}(\R\times A_R)}\lesssim R^{-\frac{n-1}{2}[1-\frac{2n}{q(n-1)}]}\|f\|_{L^p(S,d\sigma)},
\end{equation}
and the error term estimate\begin{equation}\label{eq:err-estmt}
 \|\mathcal{E}f\|_{L^q_{t,x}(\R\times A_R)}\lesssim R^{-\frac {n+1}2+\frac{n}{q}} \|f\|_{L^p(S,d\sigma)}.
\end{equation}
\end{proposition}
\begin{proof} To prove the main term estimate \eqref{eq:main-estmt}, we first observe that it is
sufficient to obtain the same estimate with the first term in the expression of $\M f$. Then
by changing to polar coordinates and the Hausdorff-Young inequality in $t$ when $q>2$ or the
Plancherel theorem in $t$ when $q=2$, we obtain
\begin{align*}
\|\M f\|_{L^q_{t,x}(\R\times A_R)}&\sim \left(\int_{R/2}^R\int_{\R}\left|r^{-\frac{n-1}{2}}\int_{I}
F(s)s^{\frac{n-3}{2}}e^{i(rs+ts)}ds\right|^qdtr^{n-1}dr\right)^{1/q}\\
&=R^{-\frac{n-1}{2}+\frac{n-1}{q}}\left(\int_{R/2}^R\int_{\R}\left|\int_{I}
F(s)s^{\frac{n-3}{2}}e^{irs}e^{its}ds\right|^qdtdr\right)^{1/q}\\
&\lesssim R^{-\frac{n-1}{2}+\frac{n-1}{q}}\left(\int_{R/2}^R\|F\|^{q}_{L^{q'}_s(I)}dr\right)^{1/q}
\lesssim R^{-\frac{n-1}{2}+\frac{n}{q}}\|f\|_{L^p(S,d\sigma)}.
\end{align*}
Hence \eqref{eq:main-estmt} follows.

To prove the error term estimate \eqref{eq:err-estmt}, for $r\ge 1$, we set
$$E(r)=\int_{0}^{\infty}e^{-ry}y^{\frac{n-3}{2}}[(y\pm2i)^{\frac{n-3}{2}}-(\pm2i)^{\frac{n-3}{2}}]dy.$$
By a similar argument as proving \cite[Chapter 3, Lemma 3.11]{Stein-Weiss:1971:fourier-analysis} (details can also be found in \cite[Proposition 3.3]{Shao:2008:paraboloid}), we have
\begin{equation}\label{eq:err-bessl}
|E(r)|\lesssim r^{-\frac{n+1}{2}}, \text{ for } r\ge 1.
\end{equation}
By changing to polar coordinates, the left-hand side of \eqref{eq:err-estmt} is comparable to
\begin{equation*}
\left(\int_{R/2}^R\int_{\R}\left|\int_{I}F(s)s^{n-2}e^{its\pm
irs}E(rs)ds\right|^qdt\,r^{n-1}dr\right)^{1/q}.
\end{equation*}
Then by the Hausdorff-Young inequality in $t$ when $q>2$ or Plancherel theorem in $t$ when $q=2$
and $s\sim 1$, it is further bounded by
\begin{equation*}
\left(\int_{R/2}^R\left|\int_{I}\left|F(s)s^{n-2}E(rs)\right|^{q'}ds\right|^{q/q'}r^{n-1}dr\right)^{1/q}
\end{equation*}
By using \eqref{eq:err-bessl} and H\"older since $q\ge p'$, it is bounded by
$R^{-\frac {n+1}2+\frac{n}{q}}\|F\|_{L^p(I)}$. Then \eqref{eq:err-estmt} follows because $\|F\|_{L^p(I)}\sim \|f\|_{L^p(S,d\sigma)}$.
\end{proof}

From the triangle inequality, we have
\begin{corollary}[Dyadic restriction estimate]\label{cor:big-R}
Suppose $f\in \L_1$. Then for all $1\le p\le \infty$, $q\ge \max\{2, p'\}$, a dyadic number $R\ge
2$ and $f\in L^p(S,d\sigma)$, we have
\begin{equation*}
\|(fd\sigma)^\vee\|_{L^q_{t,x}(\R\times A_R)}\lesssim R^{-\frac{n-1}{2}[1-\frac{2n}{q(n-1)}]}\|f\|_{L^p(S,d\sigma)}.
\end{equation*}
\end{corollary}
Having done all the preparations, we now prove Theorem \ref{thm:line-con} via the dyadic restriction estimate above.
\begin{proof}[The first proof of Theorem \ref{thm:line-con}]
We only need to show the ``sufficient" part of the claim.
We first observe that it suffices to prove \eqref{eq:lin-restr} under the boundary condition
$q>\frac{2n}{n-1}$ and $\frac{n+1}{q}=\frac{n-1}{p'}$ since other estimates are easily obtained by
a standard argument of using the H\"older inequality. From Corollary \ref{cor:big-R} and Proposition \ref{prop:sml-R}, we obtain that, for $q>\frac{2n}{n-1}$, $\frac{n+1}{q}=\frac{n-1}{p'}$, and $f\in \L_1$,
\begin{equation*}
\|(fd\sigma)^\vee\|_{L^{q}_{t,x}(\R\times A_R)}\lesssim R^{\alpha(R)}\|f\|_{L^p(S,d\sigma)},
\end{equation*}
where
$$
\alpha(R)=\begin{cases}
{-\frac{n-1}{2}[1-\frac{2n}{q(n-1)}]}, & \text{ for } R\ge 2,\\
\frac{n-1}{q}, &\text{ for } R\le 1.
\end{cases}
$$ By scaling, when $f\in \L_M$ with $M\in 2^{\Z}$, under the condition
$\frac{n+1}{q}=\frac{n-1}{p'}$,
\begin{equation*}
\|(fd\sigma)^\vee\|_{L^{q}_{t,x}(\R\times A_{R})}\lesssim (RM)^{\alpha(RM)}\|f\|_{L^p(S,d\sigma)}.
\end{equation*}
Then for general $f$, we decompose it as follows,
\begin{equation*}
f=\sum_{M:dyadic} f1_{\{(\tau, \xi):\tau=|\xi|, M\le |\xi|\le 2M\}}=\sum_{M} f_M,
\end{equation*}
where $f_M:=f1_{\{(\tau, \xi): \tau=|\xi|, M\le |\xi|\le 2M\}}$. Hence \begin{align*}
&\|(fd\sigma)^\vee\|_{L^q_{t,x}(\R\times \R^{n-1})}=\left(\sum_{R}\|(fd\sigma)^\vee\|^q_{L^q_{t,x}
(\R\times A_R)}\right)^{1/q}\\
&=\left( \sum_{R}\|\sum_{M}(f_Md\sigma)^\vee\|^q_{L^q_{t,x}(\R\times A_R)}\right)^{1/q}
\lesssim \left(\sum_{R}\left( \sum_{M}\|(f_Md\sigma)^\vee\|_{L^q_{t,x}(\R\times A_R)}\right)^q\right)^{1/q}\\
&\lesssim \left(\sum_{R}\left( \sum_{M}(RM)^{\alpha(RM)}\|f_M\|_{L^p(S,d\sigma)}\right)^q\right)^{1/q}
\lesssim \left(\sum_{M}\|f_M\|^p_{L^p(S,d\sigma)}\right)^{1/p}\sim\|f\|_{L^p(S,d\sigma)},
\end{align*}
where $R>0$ and $M>0$ are dyadic numbers; for the last line, we have used the Schur's test since $q>\frac{2n}{n-1}>p\ge 1$ and \begin{equation*}
\sup_{R>0}\sum_{M}(RM)^{\alpha(RM)}< \infty \text{ and }
\sup_{M>0}\sum_{R}(RM)^{\alpha(RM)}< \infty.
\end{equation*} Hence Theorem \ref{thm:line-con} follows.
\end{proof}

\section{Second proof of Theorem \ref{thm:line-con}}\label{sec:Nic'pf}
To begin with the second proof, we introduce the following strengthening version of the Hausdorff-Young inequality
\cite[Chapter 4, Corollary 3.16]{Stein-Weiss:1971:fourier-analysis}.
\begin{lemma}\label{le:haus-young}
If $f\in L^p(\R^n)$, $1<p\le 2$, then $\hat{f}$ belongs to $L^{p,p'}$ and
$$ \|\hat{f}\|_{L^{p',p}}\lesssim \|f\|_{L^{p}},$$
or in its dual form, for any $f\in L^{p,p'}$,
$$\|\hat{f}\|_{L^{p'}}\lesssim \|f\|_{L^{p,p'}},$$
where $L^{p,q}$ denotes the Lorentz space for $0<p<\infty$, $0<q\le \infty$, which is defined via
the equivalence that $f\in L^{p,q}$ if and only if the norm $\|f\|_{L^{p,q}}:=\left(\frac
qp\int_{0}^{\infty}\left(\lambda|\{x\in \R^n: |f(x)|>\lambda\}|^{1/p}\right)^q\frac
{d\lambda}{\lambda} \right)^{1/q}$ is finite with the usual modification weak-$L^{p}$ when
$q=\infty$.
\end{lemma}
We also introduce the following H\"older inequality in the Lorentz spaces \cite[Chapter 5, Theorem 5.3.1]{Bergh-Lofstrom:1976-interpolation-space}.
\begin{lemma}\label{le:holder}
If $0<p_1, p_2, p<\infty$ and $0<q_1, q_2, q\le \infty$ obey $\frac1p=\frac{1}{p_1}+\frac{1}{p_2}$
and $\frac1q=\frac{1}{q_1}+\frac{1}{q_2}$, then
$$\|fg\|_{L^{p,q}}\lesssim_{p_1, p_2, q_1,q_2} \|f\|_{L^{p_1,q_1}}\|g\|_{L^{p_2,q_2}},$$ whenever the right-hand side norms are finite.
\end{lemma}
Next we will present the second proof of Theorem \ref{thm:line-con}, which is inspired by
Nicola's short proof in \cite{Nicola:2008-cone-sphere-restriction} that the restriction conjecture for the sphere in $\R^n$ implies that for cone in $\R\times\R^n$.
\begin{proof}[The second proof of Theorem \ref{thm:line-con}] As
in the first proof, it is sufficient to consider $q>\frac{2n}{n-1}$ and $\frac{n+1}{q}=\frac{n-1}{p'}$. By changing to polar coordinate,
\begin{equation*}
(fd\sigma)^\vee(t,x)=c_nr^{-\frac{n-2}{2}} \int_{0}^{\infty}e^{its}F(s)s^{\frac{n-2}{2}}J_{\frac{n-2}{2}}(sr)ds.
\end{equation*}
Then by using Lemma \ref{le:haus-young} and exchanging the norms, we have
\begin{align*}
\|(fd\sigma)^\vee\|_{L^q_{t,x}}&\lesssim
\left\|\|r^{-\frac{n-2}{2}+\frac{n-1}{q}}F(s)s^{\frac n2-1}J_{\frac{n}{2}-1}(sr)\|_{L^{q',q}_s(0,\infty)}\right\|_{L^q_r(0,\infty)}\\
&\lesssim \left\|F(s)s^{\frac n2-1}\|r^{-\frac{n-2}{2}+\frac{n-1}{q}}
J_{\frac{n}{2}-1}(rs)\|_{L^q_r(0,\infty)}\right\|_{L^{q',q}_s(0,\infty)}.
\end{align*}
We observe that for each $s>0$, the integrand is bounded by
\begin{equation}\label{eq:loc-3}
\|r^{-\frac{n-2}{2}+\frac{n-1}{q}}J_{\frac{n}{2}-1}(rs)\|_{L^q_r(1/s,\infty)}+\|r^{-\frac{n-2}{2}
+\frac{n-1}{q}}J_{\frac{n}{2}-1}(rs)\|_{L^q_r(0,1/s)}.
\end{equation}
On the one hand, from the definition of the Bessel function $$J_{\frac n2-1}(r)=\frac{(r/2)^{\frac {n-2}2}}{\Gamma((n-1)/2)\Gamma(1/2)}\int_{-1}^1e^{irs}(1-s^2)^{\frac {n-3}2}ds,$$
we obtain $$J_{\frac n2-1}(r)\lesssim r^{\frac{n-2}{2}} \text{ for } n\ge 2 \text{ and }r\le 1.$$
On the other hand, from the complete expansion  of $J_m$ when $m=\frac n2-1$ and the
bound on $E(r)$ in the proof of Proposition \ref{prop:big-R}, we have
$$|J_{\frac n2-1}(r)|\lesssim r^{-1/2}+c_nr^{\frac n2-1}r^{-\frac{n+1}{2}}\lesssim r^{-1/2},
\text{ for }n\ge 2 \text{ and }r\ge 1.$$
Hence combining these two estimates on $J_{\frac n2-1}$, we obtain
$$\eqref{eq:loc-3}\lesssim s^{\frac n2-\frac nq-1},\text{ if } q>\frac{2n}{n-1}.$$
Then by the fact that $q>p$ and Lemma \ref{le:holder},
\begin{align*}
\|(fd\sigma)^\vee\|_{L^q_{t,x}}&\lesssim \|F(s)s^{\frac{n-2}{p}}s^{-\frac{n-2}{p}+n-\frac
nq-2}\|_{L^{q',q}_s}\\
&\lesssim \|F(s)s^{\frac{n-2}{p}}s^{-\frac{n-2}{p}+n-\frac nq-2}\|_{L^{q',p}_s}\\
&\lesssim \|F(s)s^{\frac{n-2}{p}}\|_{L^{p,p}}\|s^{-\frac{n-2}{p}+n-\frac nq-2}\|_{L^{\frac{1}{1/q'-1/p},\infty}}.
\end{align*}
Note the condition $\frac{n+1}{q}=\frac{n-1}{p'}$ implies that
$-\frac{n-2}{p}+n-\frac nq-2=-(\frac 1{q'}-\frac 1p)$.
Hence $$\|s^{-\frac{n-2}{p}+n-\frac
nq-2}\|_{L^{\frac{1}{1/q'-1/p},\infty}}<\infty.$$
Therefore, by the fact that $\|F(s)s^{\frac{n-2}{p}}\|_{L^{p,p}}=\|f\|_{L^p(S,d\sigma)}$, we see that
Theorem \ref{thm:line-con} follows.
\end{proof}

\bibliography{refs}

\providecommand{\bysame}{\leavevmode\hbox to3em{\hrulefill}\thinspace}
\providecommand{\MR}{\relax\ifhmode\unskip\space\fi MR }
\providecommand{\MRhref}[2]{%
  \href{http://www.ams.org/mathscinet-getitem?mr=#1}{#2}
}
\providecommand{\href}[2]{#2}
\begin{thebibliography}{10}

\bibitem{Barcelo:1985:restic-cone}
B.~Barcelo, \emph{On the restriction of the {F}ourier transform to a conical
  surface}, Trans. Amer. Math. Soc. \textbf{292} (1985), no.~1, 321--333.
  \MR{MR805965 (86k:42023)}

\bibitem{Bergh-Lofstrom:1976-interpolation-space}
J{\"o}ran Bergh and J{\"o}rgen L{\"o}fstr{\"o}m, \emph{Interpolation spaces.
  {A}n introduction}, Springer-Verlag, Berlin, 1976, Grundlehren der
  Mathematischen Wissenschaften, No. 223. \MR{MR0482275 (58 \#2349)}

\bibitem{Nicola:2008-cone-sphere-restriction}
F.~Nicola, \emph{Slicing surfaces and fourier restriction conjecture},
  arXiv:0804.3696, Proceedings of the Edinburgh Mathematical Society, to
  appear.

\bibitem{Shao:2008:paraboloid}
S.~Shao, \emph{Sharp linear and bilinear restriction estimates for paraboloids
  in the cylindrically symmetric case}, arXiv:0706.3759, Revista Matem\'atica
  Iberoamericana, to appear.

\bibitem{Stein:1979:problems-in-harmonic}
E.~M. Stein, \emph{Some problems in harmonic analysis}, Harmonic analysis in
  Euclidean spaces (Proc. Sympos. Pure Math., Williams Coll., Williamstown,
  Mass., 1978), Part 1, Proc. Sympos. Pure Math., XXXV, Part, Amer. Math. Soc.,
  Providence, R.I., 1979, pp.~3--20. \MR{MR545235 (80m:42027)}

\bibitem{Stein:1993}
\bysame, \emph{Harmonic analysis: real-variable methods, orthogonality, and
  oscillatory integrals}, Princeton Mathematical Series, vol.~43, Princeton
  University Press, Princeton, NJ, 1993, With the assistance of Timothy S.
  Murphy, Monographs in Harmonic Analysis, III. \MR{MR1232192 (95c:42002)}

\bibitem{Stein-Weiss:1971:fourier-analysis}
E.~M. Stein and G.~Weiss, \emph{Introduction to {F}ourier analysis on
  {E}uclidean spaces}, Princeton University Press, Princeton, N.J., 1971,
  Princeton Mathematical Series, No. 32. \MR{MR0304972 (46 \#4102)}

\bibitem{Strichartz:1977}
R.~S. Strichartz, \emph{Restrictions of {F}ourier transforms to quadratic
  surfaces and decay of solutions of wave equations}, Duke Math. J. \textbf{44}
  (1977), no.~3, 705--714. \MR{MR0512086 (58 \#23577)}

\bibitem{Tao:2008:recent-progress-restric}
T.~Tao, \emph{Recent progress on the restriction conjecture},
  arXiv:math/0311181.

\bibitem{Tao:1999:Bochner-Resiez-restri}
\bysame, \emph{The {B}ochner-{R}iesz conjecture implies the restriction
  conjecture}, Duke Math. J. \textbf{96} (1999), no.~2, 363--375. \MR{MR1666558
  (2000a:42023)}

\bibitem{Tao:2001:endpoint-cone}
\bysame, \emph{Endpoint bilinear restriction theorems for the cone, and some
  sharp null form estimates}, Math. Z. \textbf{238} (2001), no.~2, 215--268.
  \MR{MR1865417 (2003a:42010)}

\bibitem{Tao:2003:paraboloid-restri}
\bysame, \emph{A sharp bilinear restrictions estimate for paraboloids}, Geom.
  Funct. Anal. \textbf{13} (2003), no.~6, 1359--1384. \MR{MR2033842
  (2004m:47111)}

\bibitem{Tao-Vargas-Vega:1998:bilinear-restri-kakeya}
T.~Tao, A.~Vargas, and L.~Vega, \emph{A bilinear approach to the restriction
  and {K}akeya conjectures}, J. Amer. Math. Soc. \textbf{11} (1998), no.~4,
  967--1000. \MR{MR1625056 (99f:42026)}

\bibitem{Wolff:2001:restric-cone}
T.~Wolff, \emph{A sharp bilinear cone restriction estimate}, Ann. of Math. (2)
  \textbf{153} (2001), no.~3, 661--698. \MR{MR1836285 (2002j:42019)}

\end{thebibliography}
\bibliographystyle{amsplain}

\end{document}